\documentclass[12pt]{article}

\def\ra{\rightarrow}

\def\ss{\subseteq}

\def\Aut{{\rm Aut}}
\def\boac{{\rm boundary orbit accumulation point}}

 \def\HollowBox #1#2{{\dimen0=#1 \advance\dimen0 by -#2       
       \dimen1=#1 \advance\dimen1 by #2                       
        \vrule height #1 depth #2 width #2                    
        \vrule height 0pt depth #2 width #1                   
        \llap{\vrule height #1 depth -\dimen0 width \dimen1}%
       \hskip -#2                                             
       \vrule height #1 depth #2 width #2}}                   
 \def\BoxOpTwo{\mathord{\HollowBox{6pt}{.4pt}}\;}             

\def\endpf{\hfill $\BoxOpTwo$}

\font\teneufm=eufm10
\font\seveneufm=eufm7
\font\fiveeufm=eufm5
\newfam\eufmfam
\textfont\eufmfam=\teneufm
\scriptfont\eufmfam=\seveneufm
\scriptscriptfont\eufmfam=\fiveeufm

\newfam\msbfam
\font\tenmsb=msbm10 scaled \magstep1 \textfont\msbfam=\tenmsb
\font\sevenmsb=msbm7 scaled \magstep1 \scriptfont\msbfam=\sevenmsb
\font\fivemsb=msbm5  scaled \magstep1  \scriptscriptfont\msbfam=\fivemsb
\def\Bbb{\fam\msbfam \tenmsb}

\def\CC{{\Bbb C}}

\def\ZZ{{\Bbb Z}}

\newtheorem{theorem}{Theorem}[section]

\newtheorem{proposition}[theorem]{Proposition}

\newtheorem{example}[theorem]{EXAMPLE}

\makeindex

\usepackage{graphicx}

\usepackage{amsmath}

\begin{document}

\begin{center}
\huge \bf
Geometric Properties of Boundary Orbit Accumulation Points\footnote{{\bf MR Classification
Numbers:}  32M05, 32M17, 32T17.}
\end{center}
\vspace*{.12in}

\begin{center}
\large Steven G. Krantz\footnote{Author supported in part
by the National Science Foundation and by the Dean of the Graduate
School at Washington University.}\end{center}
\vspace*{.15in}

\begin{center}
\today
\end{center}
\vspace*{.2in}

\begin{quote}  \bf
Dedicated to Leon Ehrenpreis, a fine mathematician
and a wonderful human being.
\end{quote}
\vspace*{.2in}

\begin{quotation}
{\bf Abstract:} \sl
We study the automorphism group action on a bounded
domain in $\CC^n$.  In particular, we consider boundary
orbit accumulation points, and what geometric properties
they must have.  These properties are formulated in
the language of Levi geometry.
\end{quotation}
\vspace*{.25in}

\setcounter{section}{-1}

\section{Introduction}

In this paper a {\it domain} $\Omega \ss \CC^n$ is a connected open set.
We let ${\cal O}(\Omega)$ denote the algebra of holomorphic functions
on $\Omega$.  Also we let $\Aut(\Omega)$ be the group (under composition
of mappings) of biholomorphic self-maps of $\Omega$.  The standard
topology on $\Aut(\Omega)$ is that of uniform convergence on
compact sets (equivalently, the compact-open topology).  

We shall use the following notation:  $D$ denotes
the unit disc in the complex plane.  We let $D^2 = D \times D$ denote the bidisc,
and $D^n = D \times D \times \cdots \times D$ the polydisc in $\CC^n$.  The symbol
$B = B_n $ is the unit ball in $\CC^n$. 

Certainly domains with transitive automorphism group are of some interest.
But they are relatively few in number (see the classification theory
of Cartan, as described in [HEL]).  A very natural and compelling alternative
is to study domains with {\it noncompact automorphism group}.   A bounded
domain $\Omega$ has noncompact automorphism group if there is a sequence $\varphi_j \in \Aut(\Omega)$
such that no subsequence converges uniformly on compact sets to another automorphism.  Obversely,
the automorphism group is compact if every sequence $\{\varphi_j\}$ in $\Aut(\Omega)$
has a subsequence that converges uniformly on compact sets to another automorphism.  In this
regard, the following result of H. Cartan is central and useful (see [NAR]):

\begin{theorem} \sl
Let $\Omega \ss \CC^n$ be a bounded domain.  Then $\Omega$ has noncompact
automorphism group if and only if there are a point $P \in \partial \Omega$
and a point $X \in \Omega$ and automorphisms $\varphi_j \in \Aut(\Omega)$ such
that $\lim_{j \ra \infty} \varphi_j(X) = P$.
\end{theorem}

A point $P \in \partial \Omega$ is called a {\it boundary orbit accumulation point}
if there is a point $X \in \Omega$ and a sequence $\varphi_j \in \Aut(\Omega)$ such
that $\lim_{j \ra \infty} \varphi_j(X) = P$.   Of special interest is the
case when there is a {\it single} automorphism $\psi$ such that
$\lim_{j \ra \infty} \psi^j(X) = P$.  (Here $\psi^j$ denotes the composition
of $\psi$ with itself $j$ times when $j = 0, 1, 2, \dots$; also, if $j <0$, then
$\psi^j$ denotes the composition of $\psi^{-1}$ with itself $|j|$ times.)
In this latter
circumstance we call $P$ a {\it special boundary orbit accumulation point}.
It is not clear when an arbitrary boundary orbit accumulation point is a special
boundary orbit accumulation point.

We shall use Section 1 to collect some simple, preliminary
results that have independent interest.

\section{Background Results}

The result that inspires the present work comes from [GRK1]:

\begin{theorem} \sl
Let $\Omega \ss \CC^n$ be a smoothly bounded domain.   Suppose
that $P \in \partial \Omega$ is a boundary orbit accumulation point.
Then $P$ is a point of Levi pseudoconvexity.	   
\end{theorem}

There is an analogous result for domains without smooth boundary (and
in which the conclusion involves {\it Hartogs} pseudoconvexity).  But we
shall have no use for it in the present paper.  See [GRK1] for the
details.

We first give an example to emphasize that, even though the \boac\ is
pseudoconvex, nearby points need not be.

\begin{example}  \rm
Let $B \ss \CC^n$ be the unit ball with defining function
$\rho(z) = |z|^2 -  1$ (see [KRA1] for the concept of definining
function).   Let $\phi$ be a $C^\infty_c$ function on $\CC^n$ with
these properties:
\begin{enumerate}
\item[{\bf (a)}]  $\phi$ is real-valued and $0 \leq \phi(z) \leq 1/10$ for all $z \in \CC^n$.
\item[{\bf (b)}]  $\phi$ is radial about the point $(i,0)$.
\item[{\bf (c)}]  $\hbox{\rm supp}\, \phi \ss B((i,0), 1/10)$.
\item[{\bf (d)}]  $\phi(z) = 1/10$ for $|z - (i,0)| < 1/20$.
\end{enumerate}
Set 
$$
\Omega' = \{z \in \CC^2: - 1 + |z|^2 + \phi(z) <0\}
$$
and
$$
\Omega = \bigcap_{j=-\infty}^\infty \Phi_{1/2}^j(\Omega') \, ,
$$
where
$$
\Phi_a(z_1, z_2) = \left ( \frac{z_1 - a}{1 - \overline{a} z_1}, \frac{\sqrt{1 - |a|^2}z_2}{1 - \overline{a}z_1} \right ) \, ,
$$
any $a \in \CC$, $|a| < 1$.   It is easy to check, by direct calculation, that $\Phi_a$ is an automorphism
of the unit ball $B \ss \CC^n$.  And the domain $\Omega$ 
will be the unit ball with infinitely many strongly pseudoconcave dents that
accumulate at the points $(1,0)$ and $(-1,0)$.

Now it is plain that the point $(1,0) \in \partial \Omega$ is a \boac.
In fact we may let $X = (0,0)$ and $\varphi_j(z) = \Phi^{-j}_{1/2}(z)$ for $j = 1, 2, \dots$.  So $(1,0) \in \partial \Omega$
is certainly pseudoconvex.  Notice that, at points along the normal line through $(i,0)$, $- 1 + |z|^2 + \phi(z)$
is negative when $z$ is at least $1/5$ units from the boundary of the ball $B$.  And $-1 + |z|^2 + \phi(z)$ is positive
at $(i,0)$.  So there must be an intermediate point $\widetilde{z}$ on this line segment---a point {\it inside} the unit ball $B$---where
$-1 + |z|^2 + \phi(z)$ vanishes.  It follows that $\widetilde{z}$ is a boundary point of $\Omega'$.  Hence
$\Phi_{1/2}^j(\widetilde{z})$ is a boundary point for each $j$.
The boundary points $\Phi_{1/2}^j(\widetilde{z})$ will be strictly pseudoconcave.
\end{example}

It must be noted that, in the last example, $\Omega$ does {\it not}
have smooth boundary.  In fact, at the boundary points $(1,0)$, $(-1,0)$, the boundary
is only Lipschitz.

Of course the disc $D \ss \CC$ has noncompact automorphism group.  Let
$$
\varphi_a(\zeta) = \frac{\zeta - a}{1 - \overline{a}\zeta}
$$
for $a \in \CC$, $|a| < 1$.  Then the automorphisms 
$$
\left \{ \frac{\zeta - (1 - 1/j)}{1 - (1 - 1/j)\zeta}: j \in \ZZ \right \}
$$
are a sequence of automorphisms of $D$ that have no subsequence converging
to an automorphism.  Indeed, any subsequence either converges to the constant
function 1 or the constant function $-1$.  It is a fact---see [KRA2]---that
any domain in $\CC$ having $C^1$ boundary and noncompact automorphism group
must be conformally equivalent to the disc.   This is true without any topological
hypotheses on the domain!    In $\CC^n$, the first correct theorem of this nature---due to Bun Wong [WON]
and Rosay [ROS]---is that any $C^2$ bounded domain in $\CC^n$ with a \boac\ that
is strongly pseudoconvex must be biholomorphic to the unit ball $B$.   It is not
known in general which smoothly bounded domains have noncompact automorphism
group.   Certainly there are finite type domains with non-compact automorphism
group---see [BEP1], [BEP2].  More on this matter in what follows.

In this paper we concentrate mainly on bounded domains.  But we shall make a few
remarks right now about unbounded domains.  

\begin{example} \rm
In the complex plane $\CC$, there are
unbounded domains with noncompact automorphism group that are not the disc. 
The simplest example is when the domain $\Omega$ is the entire complex plane
$\CC$.  The punctured plane also has this property.  
\end{example}

\section{New Results}

The statement of Theorem 1.1 makes it desirable to have a formulation purely in
terms of the intrinsic, invariant geometry of the domain.  For instance, one
might hope to be able to say something about the completeness of the Kobayashi
metric at a boundary orbit accumulation point.  Unfortunately, the following
example dashes that hope:

\begin{example} \rm
Let $B \ss \CC^2$ be the unit ball with definining function
$\rho(z) = |z|^2 -  1$.  Let $\phi$ be a $C^\infty_c$ function on $\CC^n$ with
these properties:
\begin{enumerate}
\item[{\bf (a)}]  $\phi$ is real-valued and
$0 \leq \phi(z) \leq 1/10$ for all $z \in \CC^n$.
\item[{\bf (b)}]  $\phi$ is radial about the point $(i,0)$.
\item[{\bf (c)}]  $\hbox{\rm supp}\, \phi \ss B((i,0), 1/10)$.
\item[{\bf (d)}]  $\phi(z) = 1/10$ for $|z - (i,0)| < 1/20$.
\end{enumerate}
Set 
$$
\Omega' = \{z \in \CC^2: - 1 + |z|^2 + \phi(z) <0\}
$$
and
$$
\Omega = \bigcap_{-1 < a < 1}  \Phi_a (\Omega') \, ,
$$
where
$$
\Phi_a(z_1, z_2) = \left ( \frac{z_1 - a}{1 - \overline{a} z_1}, \frac{\sqrt{1 - |a|^2}z_2}{1 - \overline{a}z_1} \right ) \, ,
$$
any $a \in \CC$, $|a| < 1$.   
Then $\Omega$ is the unit ball with a groove stretching from
$(-1,0)$ to $(1,0)$.  This new domain is strongly pseudoconcave
along an entire curve from $(-1,0)$ to $(1,0)$.  Of course the point $(1,0)$ is
still a \boac.   Indeed the automorphisms $\Phi_a$, $-1 < a < 0$, send $(0,0)$ to $(1,0)$.
And, along the curve $\gamma(t) = (t,0)$, $0 < t < 1$, the Kobayashi distance to the boundary point
is infinite.

Now write $\Phi_a(z)= (\varphi_a^1(z), \varphi_a^2(z))$.
Let the point $\widetilde{z}$ be as at the end of Example 1.2.  
We consider, in the same domain $\Omega'$ as above, the curve
$$
t \longmapsto \bigl ( [1 - (1 - t)^4]\Phi_t(\widetilde{z}), 0 \bigr ) \ \ , \qquad 0 \leq t \leq 1 \, .
$$
Then this curve terminates at $(-1,0)$ and is tangent to the boundary to fourth order at
that point.  For $Q \in \partial \Omega$, let $\nu_Q$ denote the outward Euclidean unit normal
vector at $Q$.  If we take advantage of the estimate in [KRA3], to the effect that,
near a strongly pseudoconcave boundary point $Q$, a point $Q^* = Q - \delta \nu_Q$ satisfies
the estimate
$$
F_K^\Omega (Q^*, \nu) \approx \delta^{-3/4} \, ,
$$
and moreover, in the tangential direction,
$$
F_K^\Omega (Q^*, \tau) \approx C \, ,
$$
then it is easy to see that the distance to the boundary point $(-1,0)$ along the curve $\gamma$ is finite.
Therefore, at least on a domain with Lipschitz boundary, it is {\it not} the case that
a \boac\ will be a point at which the Kobayashi metric is complete.
\end{example}

On the positive side, we can prove the following result:

\begin{proposition} \sl
Let $\Omega \ss \CC^2$ be a smoothly bounded domain and $P \in \partial \Omega$
a boundary orbit accumulation point that is of finite type.  Then it is {\it not} possible for $\Omega$ to have Levi pseudoconcave boundary
points.
\end{proposition}
{\bf Proof:}  Assume that $P$ is the limit of $\varphi_j(X)$ for some point
$X \in \Omega$ and some automorphisms $\varphi_j$.  Since $P$ is of finite type
and the complex dimension is 2, we may use results of [CAT] and [ALA]
to see that curves in $\Omega$ terminating at $P$ must have infinite length
in the Kobayashi metric.

If $Q$ is a strongly pseudoconcave boundary point, obviously distinct from $P$, then we know from the
estimates in [KRA3] that the Kobayashi distance of $Q$ to $X$ is finite.
Then the $\varphi_j$ applied to $Q$ must move $Q$ into the interior
of $\Omega$ (because $X$ is being pushed infinitely far towards $P$).  
Hence $\varphi$ is {\it not} an automorphism.  That is
a contradiction.    We conclude that $Q$ cannot exist.
\endpf
\smallskip \\

A well-known conjecture in the subject says this:

\begin{quote}
{\bf The Greene-Krantz Conjecture:}  Let $\Omega$ be a smoothly bounded
domain in $\CC^n$.   If $P \in \partial \Omega$ is a \boac, then
$P$ is a point of finite type in the sense of Kohn/D'Angelo/Catlin.
\end{quote}

This conjecture has not been established in full generality.  But results
in [KIM] and [KIK1] support the conjecture.   Now we have

\begin{proposition} \sl
Let $\Omega \ss \CC^2$ be a smoothly bounded, pseudoconvex domain.  Let $P \in \partial \Omega$
be a \boac.  Assume that the Greene-Krantz conjecture is true.  Then any path
ending at $P$ will have infinite length in the Kobayashi metric.
\end{proposition}
{\bf Proof:}  This result is almost obvious.  For the hypothesis implies that
$P$ is of finite type.  And now the estimates on the Kobayashi metric in [CAT],
together with the calculations in [ALA], give the result about infinite length of
paths.
\endpf
\smallskip \\

\section{A Boundary Orbit Accumulation Point Characterization of Domains}

In [ROS] and [WON], for instance, it is shown that if a bounded domain
has a strongly pseudoconvex boundary orbit accumulation point then that
domain must be biholomorphic to the unit ball in $\CC^n$.  Put in other
words, if two distinct bounded domains have boundary orbit accumulation points,
and if those boundary orbit accumulation points are both strongly pseudoconvex,
then the two domains must be biholomorphic (since they are both biholomorphic to the ball).
 
One might more generally formulate this question
\begin{quote}
Suppose that $\Omega_1$ and $\Omega_2$ are two bounded domains in $\CC^n$.  Assume
that $\Omega_1$ has boundary orbit accumulation point $P_1$ and
$\Omega_2$ has boundary orbit accumulation point $P_2$.  If $P_1$ and
$P_2$ have the same Levi geometry, may we conclude that $\Omega_1$
is biholomorphic to $\Omega_2$?
\end{quote}

I do not know the full answer to this question at this time.  However, the
following partial answer may be proved using known techniques:

\begin{proposition} \sl
Let $\Omega_1$, $\Omega_2$ be bounded domains 
in $\CC^n$.  Let $P_1 \in \partial \Omega_1$
and $P_2 \in \partial \Omega_2$ each be boundary orbit accumulation points.
Assume that $\partial \Omega_j$ is smooth near $P_j$, $j = 1, 2$.
Suppose that each $P_j$ is of finite type in the sense of Kohn/Catlin/D'Angelo
and is also a peak point.  Finally assume that there is a neighborhood
$U_1$ of $P_1$ and a neighborhood $U_2$ of $P_2$ and a biholomorphic mapping
$$
\Phi: U_1 \cap \Omega_1 \rightarrow U_2 \cap \Omega_2
$$
such that {\bf (i)}  $\Phi$ continues to a diffeomorphism of $\partial \Omega_1 \cap U_1$
to $\partial \Omega_2 \cap U_2$, {\bf (ii)}  $\Phi(P_1) = P_2$.

Then $\Omega_1$ is biholomorphic to $\Omega_2$.
\end{proposition}
{\bf Proof:}  Choose a point $X_1 \in \Omega_1$ and automorphisms $\varphi_j$ of
$\Omega_1$ so that $\varphi_j(X_1) \rightarrow P_1$.  Likewise choose a point
$X_2 \in \Omega_2$ and automorphisms $\psi_j$ of $\Omega_2$ such that
$\psi_j(X_2) \rightarrow P_2$.	 A standard argument (see [KRA1, Chapter 11) shows that,
for any compact set $K \subseteq \Omega_1$, $\varphi_j(z)$ converges to $P_1$ 
uniformly for $z \in K$.  A similar statement holds for $\Omega_2$.

Let $K$ be a large compact set inside $\Omega_1$.  Choose $j$ so large that $\varphi_j(K) \ss U_1 \cap \Omega_1$.
Likewise let $L$ be a large compact set inside $\Omega_2$.  Choose $k$ so large that $\psi_k(L) \ss U_2 \cap \Omega_2$.
Let $\epsilon > 0$ be small and set $U_1^\epsilon = \{z \in U_1: \hbox{dist}(z, {}^c U_1 > \epsilon\}$.
Similarly set $U_2^\epsilon = \{z \in U_2: \hbox{dist}(z, {}^c U_2 > \epsilon\}$.  By shrinking $\epsilon$
if necessary, we may assume that $\Phi(\varphi_j(K)) \ss U_2^\epsilon$.  By enlarging
$L$ if necessary, we may suppose that $\varphi_k(L) \supseteq U_2^\epsilon$.

Now consider $(\psi_k)^{-1} \circ \Phi \circ \varphi_j$.  This will be a univalent holomorphic mapping
that takes $K \ss \Omega_1$ to $L$.  And the mapping is invertible.  We may similarly assume that
the inverse mapping takes $L$ to $K$.   The set of all such mappings, as $K$ exhausts $\Omega_1$ and $L$ exhausts
$\Omega_2$ forms a normal family.  And we may extract a convergent subsequence that converges to a biholomorphic
mapping of $\Omega_1$ to $\Omega_2$.  That is the result that we seek.
\endpf
\smallskip \\

There are a number of different approaches to the classical Bun Wong/Rosay theorem.   Useful
references are [WON], [ROS], [KIK2], [KIK3], [KIK4], [GKK].

\newpage

\section{Concluding Remarks}

In the paper [BEP1] and subsequent works, Bedford an Pinchuk prove
results of the following type:

\begin{quote} \sl
{\bf Theorem:}  Let $\Omega$ be a pseudoconconvex domain in $\CC^2$ with real
analytic boundary.  Suppose that $\Omega$ has noncompact automorphism
group.  Then $\Omega$ must be biholomorphic to a complex ellipsoid
of the form
$$
E_m = \{(z_1, z_2) \in \CC^2: |z_1|^2 + |z_2|^{2m} < 1\}
$$
for $m = 1, 2, \dots$.  
\end{quote}

In a private communication, David Catlin pointed out that the hypothesis
of this theorem may be reduced to ``finite type domain'' in $\CC^2$.
In later papers Bedford and Pinchuk produced analogous results in $\CC^n$.

But it must be pointed out that, in higher dimensions, we cannot hope
for a conclusion as simple as ``the domain must be a complex ellipsoid.''
For consider the domain
$$
\Omega^* = \{(z_1, z_2, z_3) \in \CC^3: |z_1|^2 + (|z_2|^2 + |z_3|^2)^2 < 1\} \, .
$$
It has automorphisms of the form
$$
\Phi_a(z_1, z_2, z_3) = \left ( \frac{z_1 - a}{1 - \overline{a} z_1}, 
\frac{\sqrt[4]{1 - |a|^2} z_2}{\sqrt{1 - \overline{a} z_1}}, \frac{\sqrt[4]{1 - |a|^2} z_3}{\sqrt{1 - \overline{a} z_1}} \right ) \, ,
$$
for $a \in \CC$, $|a| < 1$.  If we let $a$ take the values $1 - 1/j$ for $j = 1, 2, \dots$ then we see immediately
that $\Omega^*$ has noncompact automorphism group.  And $\Omega^*$ is {\it not} an ellipsoid in
the most direct sense.	 It has been conjectured by Catlin and others (see [KRA4] for the details) that
the correct conclusion in higher dimensions is that the defining function of the domain
should satisfy a certain homogeneity condition.

It has been noted that the Greene-Krantz conjecture asserts that a boundary
orbit accumulation point {\it must be} of finite type.  There is some evidence
to support the conjecture---see, for instance, [KIM], [KIK1].  If it turns
out to be true, then the Bedford/Pinchuk theorem cited above can be streamlined
to say that a smoothly bounded domain in $\CC^2$ with noncompact automorphism group
must be an ellipsiond.

It is certainly a matter of some interest to understand the
nature of boundary orbit accumulation points.	We know that
they must be pseudoconvex, and the Greene-Krantz conjecture posits
even more specific information about these points.  Another subject
of some study is boundary orbit accumulation {\it sets}---see, for instance,
[ISK] and [KRA5].  Much more can in principle be said about these sets.

Automorphism groups are in some sense an invariant that is a substitute
for the lack in several complex variables of a uniformization theorem or
a Riemann mapping theorem.  It is in our best interest to develop their
properties so that they can be used effectively to study and classify
domains up to biholomorphic equivalence.

We hope to study these matters further in future papers.
					
\newpage
	     
\null \hbox{ \ \ }

\hbox{ \ \ }

\hbox{ \ \ }

\noindent {\Large \sc References}
\bigskip  \bigskip \bigskip \\

\begin{enumerate}

\item[{\bf [ALA]}] G. Aladro, The comparability of the
Kobayashi approach region and the admissible approach region,
{\em Illinois Jour. Math.} 33(1989), 42-63.

\item[{\bf[BEP1]}] E. Bedford and S. Pinchuk, Domains in
$\CC^2$ with noncompact group of automorphisms, {\it Math.
Sb.} Nov. Ser. 135(1988), No. 2, 147-157.

\item[{\bf[BEP2]}] E. Bedford and S. Pinchuk, Domains in
$\CC^{n+1},$ with non-compact automorphism group, {\it Journal
of Geometric Analysis} 3(1991), 165--191.

\item[{\bf [CAT]}] D. Catlin, Estimates of invariant metrics on
pseudoconvex domains of dimension two, {\em Math. Z.}
200(1989), 429-466.

\item[{\bf [GKK]}] H. Gaussier, K.-T. Kim, and S. G. Krantz, A
note on the Wong-Rosay theorem in complex manifolds, {\it
Complex Var.\ Theory Appl.} 47(2002), 761--768.

\item[{\bf [GRK1]}] R. E. Greene and S. G. Krantz, Invariants
of Bergman geometry and results concerning the automorphism
groups of domains in $\CC^n$, Proceedings of the 1989
Conference in Cetraro (D. Struppa, ed.), 1991.

\item[{\bf [HEL]}] S. Helgason, {\it Differential Geometry and
Symmetric Spaces}, Academic Press, New York, 1962.

\item[{\bf [ISK]}]  A. Isaev and S. G. Krantz, On the boundary orbit accumulation 
set for a domain with non-compact automorphism group, {\it Mich.\ 
Math.\ Jour.} 43(1996), 611-617.

\item[{\bf [KIM]}]  K.-T. Kim, Domains in $\CC^n$ with a
piecewise Levi flat boundary which possess a noncompact
automorphism group, {\it Math. Ann.} 292(1992), 575--586.

\item[{\bf [KIK1]}] K.-T. Kim and S. G. Krantz, Complex scaling
and domains with non-compact automorphism group, {\it Illinois
Journal of Math.} 45(2001), 1273--1299.
		
\item[{\bf [KIK2]}] K.-T. Kim and S. G. Krantz, Some new
results on domains in complex space with non-compact
automorphism group, {\it J. Math.\ Anal.\ Appl.} 281(2003),
417--424.

\item[{\bf [KIK3]}] K.-T. Kim and S. G. Krantz, A Kobayashi
metric version of Bun Wong's theorem, {\it Complex Var.\
Elliptic Equ.} 54(2009), 355--369.

\item[{\bf [KIK4]}] K.-T. Kim and S. G. Krantz,
Characterization of the Hilbert ball by its automorphism
group, {\it Trans.\ Amer.\ Math.\ Soc.} 354(2002), 2797--2818
						 
\item[{\bf [KRA1]}] S. G. Krantz, {\it Function Theory of
Several Complex Variables}, $2^{\rm nd}$ ed., American
Mathematical Society, Providence, RI, 2001.

\item[{\bf [KRA2]}] S. G. Krantz, Characterizations of smooth
domains in $\CC$ by their biholomorphic self maps, {\it Am.
Math. Monthly} 90(1983), 555-557.

\item[{\bf [KRA3]}] S. G. Krantz, On the boundary behavior of
the Kobayashi metric, {\it Rocky Mountain J. Math.,} 22(1992).

\item[{\bf [KRA4]}]  S. G. Krantz, Recent progress and future directions in several complex
variables, in  {\it Complex Analysis Seminar}, Springer Verlag 
Lecture Notes 1268(1987), 1-23.
	
\item[{\bf [KRA5]}] S. G. Krantz, Topological/geometric
properties of the orbit accumulation set, preprint.

\item[{\bf [NAR]}] R. Narasimhan, {\it Several Complex
Variables}, University of Chicago Press, Chicago, 1971.

\newpage

\item[{\bf [ROS]}] J.-P. Rosay, Sur une characterization de la
boule parmi les domains de $\CC^n$ par son groupe
d'automorphismes, {\it Ann. Inst. Four. Grenoble} XXIX(1979),
91-97.

\item[{\bf [WON]}] B. Wong, Characterizations of the ball in
$\CC^n$ by its automorphism group, {\it Invent. Math.}
41(1977), 253-257.

\end{enumerate}
\vspace*{.95in}

\noindent \begin{quote}
Department of Mathematics \\
Washington University in St. Louis \\
St.\ Louis, Missouri 63130 \\ 			   
{\tt sk@math.wustl.edu}
\end{quote}

\end{document}